\documentclass[11pt,a4paper,reqno]{article}

\usepackage{graphicx}
\usepackage{amsmath}
\usepackage{amssymb}
\usepackage{color}

\usepackage{textcomp}
\usepackage{eurosym}
\usepackage{amsthm}
\usepackage{cite}

\usepackage[arrow,matrix,curve]{xy}

\def\figurename{Figure}
\makeatletter
\renewcommand{\fnum@figure}[1]{\textbf{\figurename~\thefigure}}

\numberwithin{equation}{section}

\newtheorem{theorem}{Theorem}[section]

\theoremstyle{definition}
\newtheorem{definition}{Definition}[section]

\theoremstyle{remark}
\newtheorem{remark}{Remark}[section]

\begin{document}

\title{\textbf{Concentric Circles Each Passing Through One Vertex of Each of Two Regular Polygons}}
\author{\textbf{Mamuka Meskhishvili}}

\date{}
\maketitle

\begin{abstract}
Given a regular $n$-gon on the plane, it is evident that from any point on the plane, taken as a center, one can draw $n$ concentric circles such that each circle passes through one of the vertices of the polygon. Naturally, this raises the problem of whether such a construction is possible for any two given regular $n$-gons on the plane. In this paper, we establish the necessary and sufficient conditions for the existence of $n$ concentric circles such that each circle passes through one vertex of each of the two regular $n$-gons.



\end{abstract}

\bigskip
\bigskip

\section{Introduction}
\label{sec:1}

Let a regular $n$-gon be given in the Euclidean plane. Denote by   \linebreak $d_1,d_2,\dots,d_n$ the distances from an arbitrary point $M$ in the plane to the vertices of the polygon, by $L$ the distance from this point to the centroid of the polygon, and by $R$ the radius of the circle circumscribed about the polygon. According to the polygonal distance theorem \cite{1, 2}, there exists a relationship among these quantities:
\begin{equation}\label{eq:1.1}
    \sum_{i=1}^n d_i^{2m}=n\bigg[(R^2+L^2)^m+\sum_{k=1}^{\lfloor\frac{m}{2}\rfloor} \binom{m}{2k}\binom{2k}{k} (RL)^{2k}(R^2+L^2)^{m-2k}\bigg],
\end{equation}
where $m=1,\dots,n-1$.

In this system, $R$ and $L$ enter symmetrically; therefore, if the ordered pair $(R,L)$ satisfies the system, then the pair $(L,R)$ is also a solution. Consequently, if $d_1,d_2,\dots,d_n$ represent the polygonal distances corresponding to a certain regular $n$-gon, then there exists another regular $n$-gon possessing the same set of polygonal distances. The radii $R_1$ and $R_2$ of the circumscribed circles of these two polygons satisfy the following system of relations \cite{5}:
\begin{equation}\label{eq:1.2}
    \sum_{i=1}^n d_i^{2m}=n\bigg[(R_1^2+R_2^2)^m+\sum_{k=1}^{\lfloor\frac{m}{2}\rfloor} \binom{m}{2k}\binom{2k}{k} (R_1R_2)^{2k}(R_1^2+R_2^2)^{m-2k}\bigg],
\end{equation}
where $m=1,\dots,n-1$.

For both polygons $P_n(R_1)$ and $P_n(R_2)$, let us take for point $M$ (the point from which the polygonal distances to the vertices are measured) the same. As shown in \cite{3, 4}:
\begin{equation}\label{eq:1.3}
    L_1=R_2, \;\; L_2=R_1.
\end{equation}

Let us draw $n$ concentric circles centered at the point $M$ with radii $d_1,d_2,\dots,d_n$. It is clear that these radii $d_i$, can be regarded as the polygonal distances from the point $M$ to the vertices of the regular polygons $P_n(R_1)$ and $P_n(R_2)$.

It is natural to pose the following two problems:
\begin{enumerate}
\item[I.] Given two regular $n$-gons, under what conditions do there exist $n$ concentric    circles such that each circle passes through one vertex of each of the two regular $n$-gons?

\item[II.] Conversely, given $n$ concentric circles, under what conditions do there exist  two regular $n$-gons whose vertices lie one by one on these circles?
\end{enumerate}

In this paper, these problems are solved.

\bigskip
\section{General Theorems}
\label{sec:2}

Condition \eqref{eq:1.3} means that the point $M$ is located at a distance $R_2$ from the center $O_1$ of the first polygon $P_1(R_1)$, and at a distance $R_1$ from the center $O_2$ of the second polygon $P_2(R_2)$. Therefore, if two regular $n$-gons are given, in order to determine the point from which the set of distances to the vertices is identical for both polygons, auxiliary circles are introduced.

\begin{definition}\label{def:1}
For two given regular $n$-gons, an auxiliary circle is defined as a circle whose center coincides with the center of one of the polygons, and whose radius is equal to the radius of the circle circumscribed about the other polygon.
\end{definition}

In the general case, the intersection of these auxiliary circles consists of two points, $M_1$ and $M_2$. Both points are equivalent; therefore, in what follows, we consider only one of them, $M_1$.

It is evident that the two given regular $n$-gons uniquely determine the radii $R_1$ and $R_2$ of their respective circumscribed circles. Furthermore, once the point $M_1$ is fixed by the intersection of the auxiliary circles, the corresponding distances $L_1$ and $L_2$    are also determined. Hence, for both polygons, all the right-hand sides of the equations in system \eqref{eq:1.1} are equal, so left-hand sides must be equal, too.

If we denote by $d_1,d_2,\dots,d_n$ and $t_1,t_2,\dots,t_n$ the polygonal distances corresponding to the two regular $n$-gons respectively, then the system \eqref{eq:1.1} yields a system of $n-1$ equations:
$$  \sum_{i=1}^n d_i^{2m}=\sum_{i=1}^n t_i^{2m},        $$
where $m=1,\dots,n-1$.

For equalization of the distances, it is sufficient that any two of the distances be equal \cite{3}. For example, without loss of generality, if we take $d_1=t_1$, then the set of all other distances will also be equal:
$$  \{d_2,\dots,d_n\}=\{t_2,\dots,t_n\}.        $$
The equality of two distances can be achieved geometrically by rotating one of the regular polygons around its own centroid \cite{4}. Therefore, we obtain:

\begin{theorem}\label{th:2.1}
Given two regular $n$-gons, there exist $n$ concentric circles, each passing through one vertex of each polygon, if and only if the following two conditions are satisfied:
\begin{enumerate}
\item[{\rm I.}] The auxiliary circles have a common point of intersection;

\item[{\rm II.}] The distances from the point of intersection to one pair of vertices, one from each polygon, are equal.
\end{enumerate}
\end{theorem}

\begin{figure}[h]
\centerline{\includegraphics[width=7cm]
    {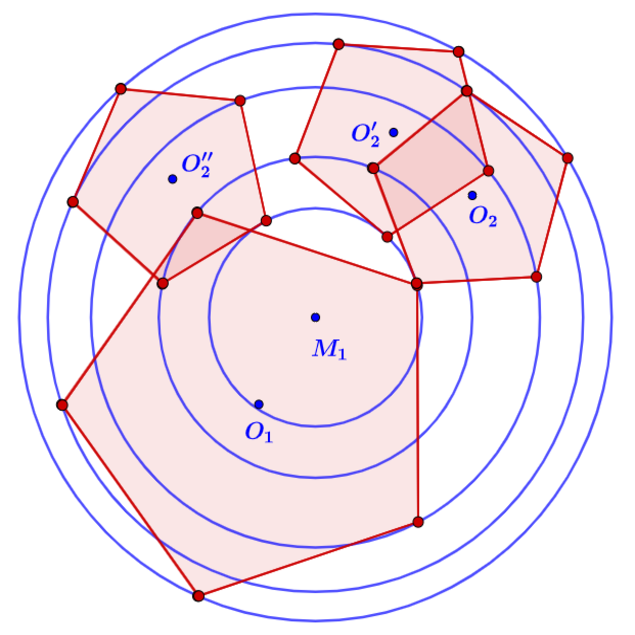} }
\caption{}
\label{fig:1}
\end{figure}

\begin{figure}[h]
\centerline{\includegraphics[width=7cm]
    {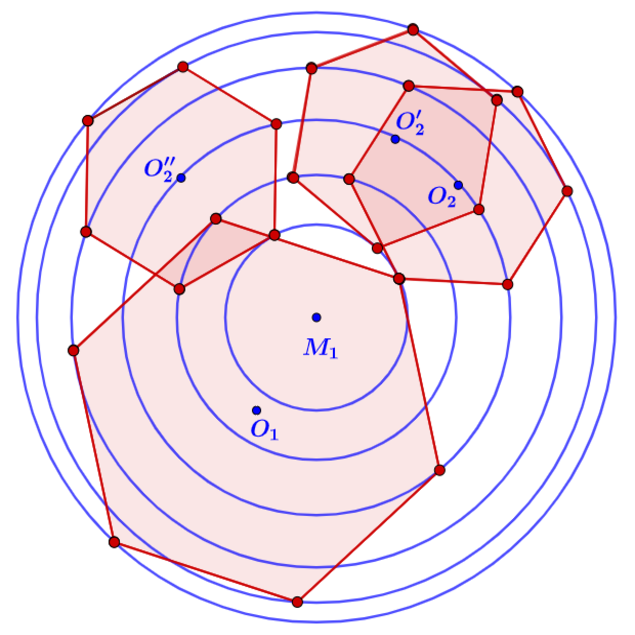} }
\caption{}
\label{fig:2}
\end{figure}

In the general case there are two sets of $n$ concentric circles corresponding to two intersection points $M_1$ and $M_2$ of auxiliary circles. Any property satisfied by one set of concentric circles is also satisfied by any other set of concentric circles with respect to the properties discussed in this article. Therefore, we will limit our discussion on a single set centered at $M_1$ (see Figures \ref{fig:1} and \ref{fig:2}).

In the case of two polygons with a shared vertex, the second condition of the Theorem  \ref{th:2.1} is automatically satisfied. As for the first condition, it is evident that the auxiliary circles have a common intersection point if the distance between the centroids of the polygons, $O_1O_2$ satisfies the following inequality:
$$  |R_1-R_2|\leq O_1O_2\leq R_1+R_2,       $$
which is automatically fulfilled in the case of a shared vertex.

Thus, we obtain the following:

\begin{theorem}\label{th:2.2}
Given two regular $n$-gons with a shared vertex, there   \linebreak  exist two sets of $n$ concentric circles such that each circle in each set passes through exactly one vertex of each of the two regular $n$-gons.
\end{theorem}

Now we consider a case when $n$ concentric circles are initially given, then their radii can be considered as polygonal distances, and the common center of the circles as the point $M_1$, i.e. in system \eqref{eq:1.2}, the left-hand side is known.

By applying cyclic averages \cite{5}:
\begin{equation}\label{eq:2.1}
    S_n^{(2m)}=\frac{1}{n}\sum_{i=1}^n d_i^{2m},
\end{equation}
where $m=1,\dots,n-1$.

Eliminating $R_1$ and $R_2$ from system\eqref{eq:1.2}, we obtain:
\begin{equation}\label{eq:2.2}
    S_n^{(2m)}=(S_n^{(2)})^m+\sum_{k=1}^{\lfloor\frac{m}{2}\rfloor} \frac{1}{2^k}\binom{m}{2k}\binom{2k}{k} \big(S_n^{(4)}-(S_n^{(2)})^2\big)^k(S_n^{(2)})^{m-2k},
\end{equation}
where $m=3,\dots,n-1$.

The first two equations of the system \eqref{eq:1.2}:
\begin{align}
    S_n^{(2)} & =R_1^2+R_2^2, \label{eq:2.3} \\
    S_n^{(4)} & =(R_1^2+R_2^2)^2+2R_1^2R_2^2, \label{eq:2.4}
\end{align}
get:
\begin{align}
    3(S_n^{(2)})^2-2S_n^{(4)} & =(R_1^2-R_2^2)^2, \label{eq:2.5} \\
    S_n^{(4)}-(S_n^{(2)})^2 & =2R_1^2R_2^2. \label{eq:2.6}
\end{align}
Equalities \eqref{eq:2.5} and (\eqref{eq:2.6} yield:
\begin{equation}\label{eq:2.7}
    \frac{2}{3}\,S_n^{(4)}\leq (S_n^{(2)})^2\leq S_n^{(4)}.
\end{equation}
Summarize \eqref{eq:2.2} and \eqref{eq:2.7} in:

\begin{theorem}\label{th:2.3}
Given $n$ concentric circles with radii $d_1,d_2,\dots,d_n$, regular $n$-gons whose vertices lie one on each of the $n$ circles exist if and only if these radii satisfy the following two conditions:
\begin{enumerate}
\medskip
\item[{\rm I.}] \quad $\displaystyle \frac{2}{3}\leq\frac{(d_1^2+d_2^2+\cdots+d_n^2)^2}{n(d_1^4+d_2^4+\cdots+d_n^4)}\leq 1$;

\bigskip
\item[{\rm II.}] \quad $\displaystyle \frac{d_1^{2m}+d_2^{2m}+\cdots+d_n^{2m}}{n}-\Big(\frac{d_1^2+d_2^2+\cdots+d_n^2}{n}\Big)^m$
\begin{multline*}
    =\sum_{k=1}^{\lfloor\frac{m}{2}\rfloor} \!\frac{1}{2^k}\binom{m}{2k}\binom{2k}{k}\Big(\frac{n(d_1^4\!+\!d_2^4\!+\cdots+\!d_n^4)\!-\!(d_1^2\!+\!d_2^2\!+\cdots+\!d_n^2)^2}{n^2}\Big)^k \\
\times \Big(\frac{d_1^2+d_2^2+\cdots+d_n^2}{n}\Big)^{m-2k},
\end{multline*}
where $m=3,\dots,n-1$.
\end{enumerate}
\end{theorem}

By solving \eqref{eq:2.5} and \eqref{eq:2.6}:

\begin{theorem}\label{th:2.4}
Given $n$ concentric circles with radii $d_1,d_2,\dots,d_n$, if there exist a regular $n$-gon whose vertices lie one on each of the $n$ circles, then also exists a second regular $n$-gon with the same property, and the circumradii of these two polygons are equal:
\begin{multline*}
    R_{1,2}^2=\frac{1}{2n}\,\Big(d_1^2+d_2^2+\cdots+d_n^2 \\
    \pm\sqrt{3(d_1^2+d_2^2+\cdots+d_n^2)^2-2n(d_1^4+d_2^4+\cdots+d_n^4)}\Big). \end{multline*}
\end{theorem}

\begin{remark}\label{rem:1}
There exists only one regular $n$-gon, if:
$$  3(d_1^2+d_2^2+\cdots+d_n^2)^2=2n(d_1^4+d_2^4+\cdots+d_n^4).   $$
\end{remark}

\bigskip
\section{Conditions for $3$ Concentric Circles and \\ Two Equilateral Triangles}
\label{sec:3}

The system \eqref{eq:1.2} yields:
\begin{align*}
    d_1^2+d_2^2+d_3^2 & =3[R_1^2+R_2^2], \\
    d_1^4+d_2^4+d_3^4 & =3\big[(R_1^2+R_2^2)^2+2R_1^2R_2^2\big].\end{align*}

The solutions for the radii of the circles are:
\begin{align*}
    d_1 & =d_1, \\
    d_2^2 & =\frac{1}{2}\,\Big(3(R_1^2+R_2^2)-d_1^2-4\sqrt{3}\,\Delta_{(R_1,R_2,d_1)}\Big), \\
    d_3^2 & =\frac{1}{2}\,\Big(3(R_1^2+R_2^2)-d_1^2+4\sqrt{3}\,\Delta_{(R_1,R_2,d_1)}\Big),
\end{align*}
where $\Delta_{(a,b,c)}$ denotes the area of a triangle with side lengths $a$, $b$, $c$.

From the first part of Theorem \ref{th:2.3} we have:
$$  (d_1^2+d_2^2+d_3^2)^2\geq 2(d_1^4+d_2^4+d_3^4),      $$
therefore
$$  \Delta_{(d_1,d_2,d_3)}\geq 0,      $$
and from  Theorem \ref{th:2.4} it follows that:

\begin{theorem}\label{th:3.1}
Let three concentric circles have radii $d_1<d_2<d_3$. Two equilateral triangles (the larger and the smaller) exist, with vertices lying one on each circle, if and only if a (nondegenerate) triangle with side lengths $d_1$, $d_2$, $d_3$ exists $(d_1+d_2>d_3)$. The radii of the circumcircles of the larger and the smaller equilateral triangles are given by:
\begin{align*}
    R_1^2 & =\frac{1}{6}\,\Big(d_1^2+d_2^2+d_3^2+4\sqrt{3}\,\Delta_{(d_1,d_2,d_3)}\Big), \\
    R_2^2 & =\frac{1}{6}\,\Big(d_1^2+d_2^2+d_3^2-4\sqrt{3}\,\Delta_{(d_1,d_2,d_3)}\Big).
\end{align*}
\end{theorem}

\begin{remark}\label{rem:2}
In the degenerate case $(d_1+d_2=d_3)$, only one equilateral triangle exists.
\end{remark}

\begin{figure}[h]
\centerline{\includegraphics[width=7cm]
    {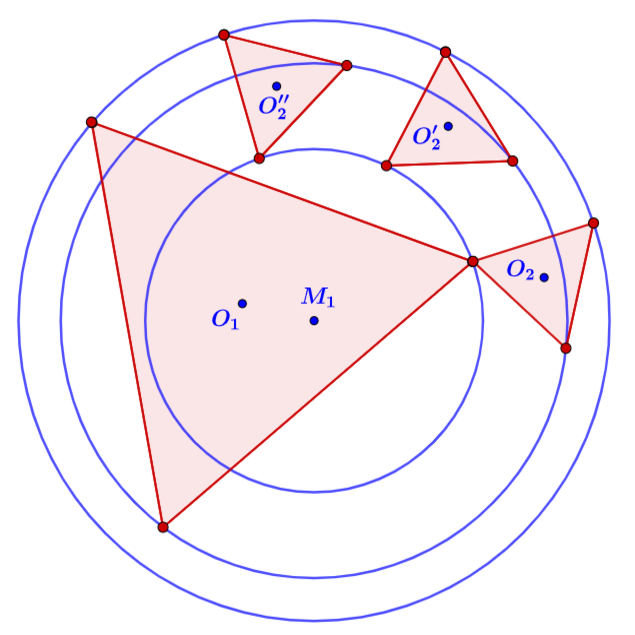} }
\caption{}
\label{fig:3}
\end{figure}

Figure \ref{fig:3} shows a large triangle with centroid $O_1$ and a small triangle ``rotating'' around it, shown in three positions (with centroids $O_2$, $O_2'$, $O_2''$). Each of the small triangle's three vertices lies on one of three concentric circles whose common center is the point $M_1$.

\bigskip
\section{Conditions for $4$  Concentric Circles and \\ Two Squares}
\label{sec:4}

The radii of the circles $d_1$, $d_2$, $d_3$, $d_4$ and the radii of circumscribed circles of two squares $R_1$, $R_2$ satisfy:
\begin{align*}
    d_1^2+d_2^2+d_3^2+d_4^2 & =4[R_1^2+R_2^2], \\
    d_1^4+d_2^4+d_3^4+d_4^4 & =4\big[(R_1^2+R_2^2)^2+2R_1^2R_2^2\big], \\
    d_1^6+d_2^6+d_3^6+d_4^6 & =4\big[(R_1^2+R_2^2)^3+6R_1^2R_2^2(R_1^2+R_2^2)\big].
\end{align*}
From the second part of the Theorem \ref{th:2.3} follows:
\begin{multline*} 8(d_1^6+d_2^6+d_3^6+d_4^6)+(d_1^2+d_2^2+d_3^2+d_4^2)^3 \\
=6(d_1^2+d_2^2+d_3^2+d_4^2)(d_1^4+d_2^4+d_3^4+d_4^4),
\end{multline*}
which is equivalent to
$$  3(d_1^2+d_2^2-d_3^2-d_4^2)(d_1^2+d_3^2-d_2^2-d_4^2)(d_1^2+d_4^2-d_2^2-d_3^2)=0. $$
Without loss of generality, let us assume that the radii of the circles satisfy:
$$  d_1<d_2<d_3<d_4,   $$
then holds:
$$  d_1^2+d_4^2=d_2^2+d_3^2.        $$

Let us use the term outer circle radii for $d_1$ and $d_4$, and inner circle radii for $d_2$ and $d_3$.

From the first part of the Theorem \ref{th:2.3}:
$$  3(d_1^2+d_2^2+d_3^2+d_4^2)^2\geq 8(d_1^4+d_2^4+d_3^4+d_4^4),        $$
on the other hand,
\begin{multline*}
    3(d_1^2+d_2^2+d_3^2+d_4^2)^2-8(d_1^4+d_2^4+d_3^4+d_4^4) \\
    =64\Delta_{(d_1,d_4,\sqrt{2}\,d_2)}^2=64\Delta_{(d_1,d_4,\sqrt{2}\,d_3)}^2=64\Delta_{(d_2,d_3,\sqrt{2}\,d_4)}^2=64\Delta_{(d_2,d_3,\sqrt{2}\,d_1)}^2.
\end{multline*}
This means that if one of the triangles with sides --
$$  (d_1,d_4,\sqrt{2}\,d_2),\;(d_1,d_4,\sqrt{2}\,d_3),\;(d_2,d_3,\sqrt{2}\,d_4),\;(d_2,d_3,\sqrt{2}\,d_1) $$
exists, then all of them exist.

Let us call them the associated triangles for the four circles. The sides of the associated triangles are composed as follows:
\begin{enumerate}
\item[{\rm (a)}] two outer circle radii and $\sqrt{2}$ times one inner circle radius;

\item[{\rm (b)}] two inner circle radii and $\sqrt{2}$ times one outer circle radius.
\end{enumerate}

The obtained results are summarized in:

\begin{theorem}\label{th:4.1}
Let four concentric circles have radii $d_1<d_2<d_3<d_4$. Two squares (the larger and the smaller) exist, with vertices lying one on each circle, if and only if the following two conditions are satisfied:
\begin{enumerate}
\item[{\rm I.}] An (nondegenerate) associated triangle exists;

\item[{\rm II.}] The sum of the squares of the outer and inner radii are equal:
$$  d_1^2+d_4^2=d_2^2+d_3^2;        $$
\end{enumerate}
The radii of the circumcircles of the larger and the smaller squares are given by:
\begin{align*}
    R_1^2 & =\frac{1}{4}\,(d_1^2+d_4^2)+\Delta_{(d_1,d_4,\sqrt{2}\,d_2)}, \\
    R_2^2 & =\frac{1}{4}\,(d_1^2+d_4^2)-\Delta_{(d_1,d_4,\sqrt{2}\,d_2)}.
\end{align*}
\end{theorem}

\begin{remark}\label{rem:3}
If the associated triangle is degenerate, only one square exists.
\end{remark}

\begin{figure}[h]
\centerline{\includegraphics[width=7cm]
    {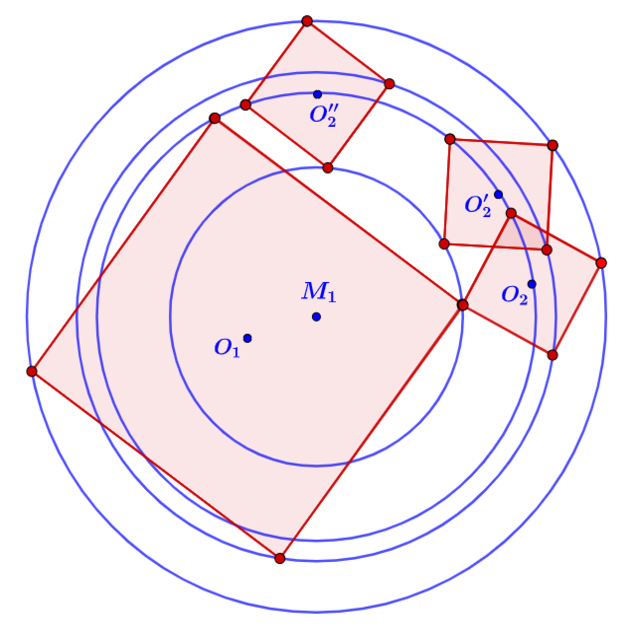} }
\caption{}
\label{fig:4}
\end{figure}

Figure \ref{fig:4} shows a large square with centroid $O_1$ and a small square   \linebreak ``rotating'' around it, shown in three positions (with centroids $O_2$, $O_2'$, $O_2''$).

Follows from the system \eqref{eq:1.2}:
\begin{align*}
    d_1^2+d_4^2=d_2^2+d_3^2 & =2(R_1^2+R_2^2), \\
    d_1^2d_4^2+d_2^2d_3^2 & =2(R_1^4+R_2^4).
\end{align*}

From which we obtain expressions for the radii of the concentric circles:
\begin{align*}
    d_1 & =d_1, \\
    d_2^2 & =R_1^2+R_2^2-4\Delta_{(R_1,R_2,d_1)}, \\
    d_3^2 & =R_1^2+R_2^2+4\Delta_{(R_1,R_2,d_1)}, \\
    d_4^2 & =2(R_1^2+R_2^2)-d_1^2.
\end{align*}

\vskip+1cm



\noindent DEPARTMENT OF MATHEMATICS

\noindent GEORGIAN-AMERICAN HIGH SCHOOL

\noindent 18 CHKONDIDELI STR., TBILISI 0180, GEORGIA

\noindent \textit{E-mail address:} \texttt{mathmamuka@gmail.com}


\begin{thebibliography}{5}


\bibitem{1} Meskhishvili, M., \textit{Cyclic averages of regular polygons and platonic solids}, Commun. Math. Appl., \textbf{11 (3) (2020)}, 335--355.

\bibitem{2} Meskhishvili, M., \textit{Cyclic averages of regular polygonal distances}, Int. J. Geom., \textbf{10 (1) (2021)}, 58--65.

\bibitem{3} Meskhishvili, M., \textit{Two regular polygons with a shared vertex}, Commun. Math. Appl., \textbf{13 (2) (2022)}, 435--447.

\bibitem{4} Meskhishvili, M., \textit{Two non-congruent regular polygons having vertices at the same distances from the point}, Int. J. Geom., \textbf{12 (1) (2023)}, 35--45.

\bibitem{5} Meskhishvili, M., \textit{Polygonal distances theorems for two regular polygons}, Int. J. Geom., \textbf{13 (4) (2024)}, 26--36.

\end{thebibliography}
\end{document}